\documentclass[12pt,
twoside,english]{amsproc}

\PassOptionsToPackage{pdfauthor={Vladimir Kisil},%
    pdftitle={Spectrum as the Support of Functional Calculus},%
    pdfsubject={Functional Analysis},%
    backref=page,%
  pdfkeywords={symmetries, functional calculus, spectrum}
}{hyperref}

\let\oldtheenumi=\theenumi
\renewcommand{\theenumi}{\textup{\oldtheenumi}}
\newtheorem{thm}{Theorem}[section]
\newtheorem{prop}[thm]{Proposition}

\theoremstyle{definition}%
\newtheorem{defn}[thm]{Definition}

\newtheorem{example}[thm]{Example}

\theoremstyle{remark}

\usepackage[noshort,extnum,printedin]{myart2k}

\ppnum{LEEDS--MATH--PURE--200\reflectbox{2}-21}{\href{http://arXiv.org/abs/math/0208249}{math/0208249}}
{200\reflectbox{2}}

\date{August 30, 200\reflectbox{2}}
\input{mydef}
\providecommand{\matr}[4]{{\ensuremath{ \left(\!\! \begin{array}{cc}
#1 & #2 \\ #3 & #4
\end{array}\!\!\right) }}}

\newcommand{\Ba}{\bar{\alpha}}
\newcommand{\Bb}{\bar{\beta}}

\providecommand{\oper}[1]{\mathcal{#1}}
\newcommand{\ga}{\mathsf{a}}

\ifundefined{eoe}
   \DeclareMathSymbol{\eoe}{\mathord}{AMSa}{"06}
\fi

\newcounter{myenumi}
\renewcommand{\themyenumi}{(\alph{myenumi})}

\providecommand{\SL}{\ensuremath{SL(2,\Space{R}{})}}

\usepackage{graphicx}

\begin{document}

\author[Vladimir V. Kisil]%
{\href{http://maths.leeds.ac.uk/~kisilv/}{Vladimir V. Kisil}}
\thanks{On  leave from Odessa University.}

\address{%            
%Institute of Mathematics\\
%Economics and Mechanics\\
%Odessa State University\\
%ul. Petra Velikogo, 2\\
%Odessa-57, 270057, UKRAINE
School of Mathematics\\
University of Leeds\\
Leeds LS2\,9JT\\
UK
}

\email{\href{mailto:kisilv@maths.leeds.ac.uk}{kisilv@maths.leeds.ac.uk}}

\urladdr{\href{http://maths.leeds.ac.uk/~kisilv/}%
{http://maths.leeds.ac.uk/\~{}kisilv/}}

\title[Spectrum as the Support of Functional Calculus]{Spectrum as the
  Support\\
 of Functional Calculus}

\begin{abstract}
  We investigate the new definition of analytic functional calculus   in the terms
  of representation theory of \(\SL\). We avoid any usage of its algebraic
  homomorphism property and replace it by the demand
  to be an intertwining operator. The related notion of spectrum and
  spectral mapping theorem are given. The construction is illustrated
  by a simple example of calculus and spectrum of non-normal \(n\times
  n\) matrix. 
\end{abstract}
%  \AMSMSC{47A60}{46H30}
\maketitle

\tableofcontents\vspace{3mm}

\setcounter{section}{-1}

\section{Introduction}
\label{sec:introduction-1}

United in the trinity functional calculus, spectrum, and spectral
mapping theorem play the exceptional r\^ole in functional
analysis and could not be substituted by anything else. 
All traditional  definitions of functional calculus  are covered by the
following rigid template based on \emph{algebra homomorphism} property:
\begin{defn}
\label{de:calculus-old}
 An \emph{functional calculus} for an element
 \(a\in\algebra{A}\)\ is a continuous 
linear mapping
\(\Phi: \mathcal{ A}\rightarrow \algebra{A}\)\ such that
\begin{enumerate} 
\item 
 \(\Phi\)\ is a unital \emph{algebra homomorphism}
 \begin{displaymath}
   \Phi(f \cdot g)=\Phi(f) \cdot \Phi (g).
\end{displaymath}
\item 
 There is an initialisation condition: \(\Phi[v_0]=a\)\ for
 for a fixed function \(v_0\), e.g. \(v_0(z)=z\). 
\end{enumerate}
\end{defn}

Most typical definition of the spectrum is seemingly independent and 
uses the important notion of
resolvent: 
\begin{defn}
  \label{de:spectrum}
  A \emph{resolvent} of element \(a\in\algebra{A}\)\ is the function
  \(R(\lambda)=(a-\lambda e)^{-1}\), which is the image under
  \(\Phi\)\ of the Cauchy kernel \((z-\lambda)^{-1}\).

  A  \emph{spectrum} of \(a\in\algebra{A}\)\ is the set \(\spec a\)\ of
  singular points of its resolvent \(R(\lambda)\).
\end{defn}
Then the following important theorem links spectrum and functional calculus
together. 
\begin{thm}[Spectral Mapping]
  \label{th:spectral-mapping}
  For a function \(f\) suitable for the   functional calculus:
   \begin{equation}\label{eq:spectral-mapping}
        f(\spec a)=\spec  f(a).
   \end{equation}
\end{thm}

However the power of the classic spectral theory rapidly decreases if
we move beyond the study of one normal operator (e.g. for
quasinilpotent ones) and is virtually nil if we consider several
non-commuting ones.
Sometimes these severe limitations are seen to be irresistible and
alternative constructions, i.e. model theory~\cite{Nikolskii86}, were
developed.

Yet the spectral theory can be revived from a fresh start. While three
compon\-ents---functional calculus, spectrum, and spectral mapping
theorem---are highly interdependent in various ways 
we will nevertheless arrange them as follows: 

\begin{enumerate}
\item Functional  calculus is an \emph{original} notion defined in
  some independent terms;
\item Spectrum (or spectral decomposition) is derived from previously
  defined functional calculus as its \emph{support} (in some
  appropriate sense);
\item Spectral mapping theorem then should drop out naturally in the
  form~\eqref{eq:spectral-mapping} or some its variation.
\end{enumerate}

Thus the entire scheme depends from the notion of the functional
calculus and our ability to escape limitations of
Definition~\ref{de:calculus-old}.  The first known to the present
author definition of functional calculus not linked to algebra
homomorphism property was the Weyl functional calculus defined by an
integral formula~\cite{Anderson69}. Then its intertwining property
with affine transformations of Euclidean space was proved as a
theorem. However it seems to be the only ``non-homomorphism'' calculus
for decades.

The different approach to whole range of calculi was given
in~\cite{Kisil95i} and developed in~\cite{Kisil98a} in terms of
\emph{intertwining operators} for group representations. It was
initially targeted for several non-commuting operators because no
non-trivial algebra homomorphism with a commutative algebra of
function is possible in this case.  However it emerged later that the
new definition is a useful replacement for classical one across all
range of problems.

In the present note we will support the last claim by consideration of
the simple known problem: characterisation a \(n \times n\)\ matrix up
to similarity. Even that ``freshman'' question could be only sorted
out by the classical spectral theory for a small set of diagonalisable
matrices. Our solution in terms of new spectrum will be full and thus
unavoidably coincides with one given by the Jordan normal form of
matrix. Other more difficult questions are the subject of ongoing
research.

\section{Another Approach to Analytic Functional Calculus} 

Anything called ``\emph{functional} calculus'' uses properties of
\emph{functions} to model properties of \emph{operators}. Thus
changing our viewpoint on functions we could get another approach to
operators. We start from the following observation reflected in the
almost any textbook on complex analysis:
\begin{prop}
\label{pr:function-theory}
\emph{Analytic function theory} in the unit disk \(\Space{D}{}\) is
a manifestation of %intimately connected with 
the \textit{mock discrete series}
representation \(\rho_1\) of \(\SL\):
\begin{equation}\label{eq:rho-1}
  \rho_1(g): f(z) \mapsto
  \frac{1}{\alpha-{\beta}{z}} \,
  f\left(
    \frac{\Ba z - \Bb}{\alpha-{\beta} z} 
  \right), \quad \textup{ where }
  \matr{\Ba}{-\Bb}{-\beta}{\alpha}\in\SL.
\end{equation}
\end{prop}
The representation~\eqref{eq:rho-1} is unitary irreducible when acts
on the Hardy space \(\FSpace{H}{2}\). Consequently we have one more
reason to abolish the template definition~\ref{de:calculus-old}:
\(\FSpace{H}{2}\) is \emph{not} an algebra. Instead we replace the
\textit{homomorphism  property} by a \textit{symmetric covariance}:
\begin{defn}
  \label{de:functional-calculus-new}
  An \emph{analytic functional calculus} for an element
  \(a\in\algebra{A}\)\ and  an
  \(\algebra{A}\)-module \(M\)\ is a \textit{continuous 
  linear} mapping
  \(\Phi:\FSpace{A}{}(\Space{D}{})\rightarrow \FSpace{A}{}(\Space{D}{},M)\)\ such that 
  \begin{enumerate} 
  \item \(\Phi\)\ is an \emph{intertwining operator} 
    \begin{displaymath}
      \Phi\rho_1=\rho_a\Phi
    \end{displaymath}
    between two representations of the
    \(\SL\)\ group \(\rho_1\)~\eqref{eq:rho-1} and \(\rho_a\)\
    defined bellow in~\eqref{eq:rho-a}.
  \item There is an initialisation condition: \(\Phi[v_0]=m\)\ for
    \(v_0(z)\equiv 1\) and \(m\in M\), where \(M\) is a left
    \(\algebra{A}\)-module.  
  \end{enumerate}
\end{defn} Note that our functional calculus released form the
homomorphism condition can take value in any left
\(\algebra{A}\)-module \(M\), which however could be \(\algebra{A}\)
itself if suitable. This add much flexibility to our construction.

The earliest functional calculus, which is \emph{not} an algebraic
homomorphism, was the Weyl functional calculus and
was defined just by an integral formula as an operator valued
distribution~\cite{Anderson69}. In that paper
(joint) spectrum was defined as support of the Weyl calculus, i.e. as
the set of point where this operator valued distribution does not
vanish. We also define
the spectrum as a support of functional calculus, but due to our
Definition~\ref{de:functional-calculus-new} it will means the set of
non-vanishing intertwining operators with primary subrepresentations.
\begin{defn}
  \label{de:spectrum-new}
    A corresponding \emph{spectrum} of \(a\in\algebra{A}\) is the
  \textit{support} of the functional
  calculus \(\Phi\), i.e. the collection of intertwining operators of
  \(\rho_a\) with \emph{prime representations}~\cite[\S~8.3]{Kirillov76}.
\end{defn}

More variations of functional calculi are obtained from other groups and their
representations~\cite{Kisil95i,Kisil98a}. 

\section{Background in Complex Analysis from $\SL$ Group}
\label{sec:backgr-compl-analys}
To understand the functional calculus from
Definition~\ref{de:functional-calculus-new} we need first to
realise the function theory from Proposition~\ref{pr:function-theory},
see~\cite{Kisil97c,Kisil97a,Kisil01a,Kisil02c} for more details.

Elements of \(\SL\) could be represented by \(2\times 2\)-matrices
with complex entries such that:
\begin{displaymath}%[eq:sl2r]
  g= \matr{\alpha}{\Bb}{-\beta}{\Ba},
  \qquad 
  g^{-1}= \matr{\Ba}{-\Bb}{-\beta}{\alpha}, 
  \qquad
  \modulus{\alpha}^2-\modulus{\beta}^2=1.
\end{displaymath}
There are other realisations of \(\SL\) which may be more suitable under
other circumstances, e.g. in the upper half-plane.

We may identify the unit disk \(\Space{D}{}\) with the left coset
\(\Space{T}{} \backslash \SL\) for the unit circle \(\Space{T}{}\) through
the important decomposition \(\SL\sim \Space{T}{}\times\Space{D}{}\)\
with \(K=\Space{T}{}\)---the only compact subgroup of \(\SL\):
\begin{equation}\label{eq:sl2-u-psi-coord}
  \matr{\alpha}{\bar{\beta}}{\beta}{\bar{\alpha}} 
%  & =& \modulus{\alpha} 
%  \matr{ 
%    \frac{{\alpha}}{ \modulus{\alpha} } }{0}{0}{\frac{\bar{\alpha}}{ 
%      \modulus{\alpha} } }
%  \matr{1}{\bar{\beta}{\alpha}^{-1}}{{\beta}\bar{\alpha}^{-1}}{1}
%  \nonumber \\
  = \frac{1}{\sqrt{1- \modulus{u}^2 }}
  \matr{e^{i\omega}}{0}{0}{e^{-i\omega}}
  \matr{1}{u}{\bar{u}}{1}, 
\end{equation}
where
\begin{displaymath}
  \omega=\arg \alpha,\qquad 
  u=\bar{\beta}{\alpha}^{-1},\qquad \modulus{u}<1.
\end{displaymath}
Each element \(g\in\SL\) acts by the linear-fractional transformation
(the M\"obius map) on \(\Space{D}{}\)\ 
and \(\Space{T}{}\)
\(\FSpace{H}{2}(\Space{T}{})\) as follows: 
\begin{equation}\label{eq:moebius}
  g^{-1}: z \mapsto \frac{\Ba z - \Bb}{\alpha-{\beta} z},
\qquad \textrm{ where } \quad
g^{-1}=\matr{\Ba}{-\Bb}{-\beta}{\alpha}.
\end{equation}
In the decomposition~\eqref{eq:sl2-u-psi-coord} the first matrix on
the right hand side acts by transformation~\eqref{eq:moebius} as an
orthogonal rotation of \(\Space{T}{}\) and \(\Space{D}{}\); and the
second one---by transitive family of maps of the unit disk onto
itself.

The standard linearisation procedure~\cite[\S~7.1]{Kirillov76} leads
from M\"obius transformations~\eqref{eq:moebius} to the unitary
representation \(\rho_1\) irreducible on the \emph{Hardy space}:
\begin{equation}\label{eq:rho-1-1}
  \rho_1(g): f(z) \mapsto
  \frac{1}{\alpha-{\beta}{z}} \,
  f\left(
    \frac{\Ba z - \Bb}{\alpha-{\beta} z} 
  \right)
\qquad  \textrm{ where } \quad 
g^{-1}=\matr{\Ba}{-\Bb}{-\beta}{\alpha}.
\end{equation}
M\"obius transformations provide a natural family of
intertwining operators for \(\rho_1\) coming from inner
automorphisms of \(\SL\) (will be used later). 
 
We choose~\cite{Kisil98a,Kisil01a} \(K\)-invariant function \(v_0(z)\equiv 1\) 
to be  a \emph{vacuum vector}.
Thus the associated \emph{coherent states}
\begin{displaymath}
  v(g,z)=\rho_1(g)v_0(z)= (u-z)^{-1}
\end{displaymath} are completely determined by the point on the unit disk \(
u=\bar{\beta}{\alpha}^{-1}\). The family of coherent states considered
as a function of both \(u\) and \(z\) is obviously the \textit{Cauchy
  kernel}~\cite{Kisil97c}. The \emph{wavelet transform}~\cite{Kisil97c,Kisil98a}
\(\oper{W}:\FSpace{L}{2}(\Space{T}{})\rightarrow
\FSpace{H}{2}(\Space{D}{}): f(z)\mapsto
\oper{W}f(g)=\scalar{f}{v_g}\)\ is the \textit{Cauchy integral}:
\begin{equation}\label{eq:cauchy}
  \oper{W} f(u)=\frac{1}{2\pi i}\int_{\Space{T}{}}f(z)\frac{1}{u-z}\,dz.
\end{equation}

Other classical objects of complex analysis (the Cauchy-Riemann
equation, the Taylor series, the Bergman space, etc.)  can be also
obtained~\cite{Kisil97c,Kisil01a} from representation \(\rho_1\) but
are not used and considered here.

\section{Representations of $\SL$ in Banach Algebras}
A simple but important observation is that the M\"obius
transformations~\eqref{eq:moebius} can be easily extended to any 
Banach algebra.
  Let \(\algebra{A}\) be a Banach algebra with the unit \(e\), 
  an element \(a\in\algebra{A}\) with \(\norm{a}<1\) be fixed, then 
  \begin{equation}\label{eq:sl2-on-A}
    g: a \mapsto g\cdot a=(\Ba a -\Bb e)(\alpha e-\beta a)^{-1}, \qquad
    g\in\SL
  \end{equation}
  is a well defined \(\SL\) action on a subset \(\Space{A}{}=\{g\cdot
  a \such g\in 
  \SL\}\subset\algebra{A}\), i.e. \(\Space{A}{}\) is a \(\SL\)-homogeneous
  space. Let us define the \emph{resolvent} function
  \(R(g,a):\Space{A}{}\rightarrow \algebra{A}\):
  \begin{displaymath}
    R(g, a)=(\alpha e-\beta a)^{-1} \quad 
  \end{displaymath}
  then 
  \begin{equation}\label{eq:ind-rep-multipl}
    R_1(g_1,\ga)R_1(g_2,g_1^{-1}\ga)=R_1(g_1g_2,\ga).
  \end{equation}
  The last identity is well known in representation
  theory~\cite[\S~13.2(10)]{Kirillov76} and is a key ingredient of
  \emph{induced representations}. Thus we can again
  linearise~\eqref{eq:sl2-on-A} (cf.~\eqref{eq:rho-1-1}) in
  the space of continuous functions \(\FSpace{C}{}(\Space{A}{},M)\)
  with values in  a left
  \(\algebra{A}\)-module \(M\), e.g.\(M=\algebra{A}\):
  \begin{eqnarray}
    \rho_a(g_1): f(g^{-1}\cdot a ) &\mapsto&
    R(g_1^{-1}g^{-1}, a)f(g_1^{-1}g^{-1}\cdot a) \label{eq:rho-a}\\
    &&\quad =
    (\alpha' e-\beta'a)^{-1} \,
    f\left(
      \frac{\Ba' \cdot a - \Bb' e}{\alpha'  e -\beta' a} 
    \right).  \nonumber
  \end{eqnarray}  
  For any \(m\in M\) we can again define a \(K\)-invariant
  \emph{vacuum vector} as \(v_m(g^{-1}\cdot 
  a)=m\otimes v_0(g^{-1}\cdot a) \in \FSpace{C}{}(\Space{A}{},M)\). 
  It generates the associated with \(v_m\) family of \emph{coherent
    states} \(v_m(u,a)=(ue-a)^{-1}m\), where \(u\in\Space{D}{}\).

The \emph{wavelet transform}  defined by
the same common formula based on coherent states (cf.~\eqref{eq:cauchy}):
\[\oper{W}_m f(g)= \scalar{f}{\rho_a(g) v_m},\qquad \]
is a version of Cauchy integral, which maps
\(\FSpace{L}{2}(\Space{A}{})\) to \(\FSpace{C}{}(\SL,M)\). It is
 closely related (but not identical!) to the
Riesz-Dunford functional calculus:  the traditional functional
calculus is given by the case:
\begin{displaymath}
  \Phi: f \mapsto \oper{W}_m f(0) \qquad\textrm{ for } M=\algebra{A}
  \textrm{ and } m=e.
\end{displaymath}

The both conditions---the intertwining property and initial
value---required by Definition~\ref{de:functional-calculus-new} easily
follows from our construction.

\section{Jet Bundles and Prolongations of $\rho_1$}
\label{sec:jet-bundl-prol-1}
Spectrum was defined in~\ref{de:spectrum-new} as
the \emph{support} of our functional calculus. To elaborate its meaning we
need the notion of a \emph{prolongation} of representations introduced by
S.~Lie, see  \cite{Olver93,Olver95} for a detailed exposition.

\begin{defn} \textup{\cite[Chap.~4]{Olver95}}
  Two holomorphic functions have \(n\)th \emph{order contact} in a point
  if their value and their first \(n\) derivatives agree at that point,
  in other words their Taylor expansions are the same in first \(n+1\)
  terms. 

  A point \((z,u^{(n)})=(z,u,u_1,\ldots,u_n)\) of the \emph{jet space}
  \(\Space{J}{n}\sim\Space{D}{}\times\Space{C}{n}\) is the equivalence
  class of holomorphic functions having \(n\)th contact at the point \(z\)
  with the polynomial:
  \begin{equation}\label{eq:Taylor-polynom}
    p_n(w)=u_n\frac{(w-z)^n}{n!}+\cdots+u_1\frac{(w-z)}{1!}+u.
  \end{equation}
\end{defn}

For a fixed \(n\) each holomorphic function
\(f:\Space{D}{}\rightarrow\Space{C}{}\) has \(n\)th \emph{prolongation}
(or \emph{\(n\)-jet}) \(\object[_n]{j}f: \Space{D}{} \rightarrow
\Space{C}{n+1}\): 
\begin{equation}\label{eq:n-jet}
  \object[_n]{j}f(z)=(f(z),f'(z),\ldots,f^{(n)}(z)).
\end{equation}The graph \(\Gamma^{(n)}_f\) of \(\object[_n]{j}f\) is a
submanifold of \(\Space{J}{n}\) which is section of the \emph{jet
bundle} over \(\Space{D}{}\) with a fibre \(\Space{C}{n+1}\). We also
introduce a notation \(J_n\) for the map \(
  J_n:f\mapsto\Gamma^{(n)}_f
\) of a holomorphic \(f\) to the graph \(\Gamma^{(n)}_f\) of its \(n\)-jet
\(\object[_n]{j}f(z)\)~\eqref{eq:n-jet}.

One can prolong any map of functions \(\psi: f(z)\mapsto [\psi f](z)\) to
a map \(\psi^{(n)}\) of \(n\)-jets by the formula
\begin{equation}\label{eq:prolong-def}
  \psi^{(n)} (J_n f) = J_n(\psi f).
\end{equation} For example such a prolongation \(\rho_1^{(n)}\) of the
representation \(\rho_1\) of the group \(\SL\) in
\(\FSpace{H}{2}(\Space{D}{})\) (as any other representation of a Lie
group~\cite{Olver95}) will be again a representation of
\(\SL\). Equivalently we can say that \(J_n\) \emph{intertwines} \(\rho_1\) and
\(\rho^{(n)}_1\):
\begin{displaymath}
   J_n \rho_1(g)= \rho_1^{(n)}(g) J_n \quad
  \textrm{ for all } g\in\SL.
\end{displaymath}
Of course, the representation \(\rho^{(n)}_1\) is not irreducible: any jet
subspace \(\Space{J}{k}\), \(0\leq k \leq n\) is
\(\rho^{(n)}_1\)-invariant subspace of \(\Space{J}{n}\).  However the
representations \(\rho^{(n)}_1\) are
\emph{primary}~\cite[\S~8.3]{Kirillov76} in the sense that they are not 
sums of two subrepresentations.

The following statement explains why jet spaces appeared in our study
of functional calculus.
\begin{prop}
  \label{pr:Jordan-zero}
  Let matrix \(a\) be a Jordan block of a length \(k\) with the
  eigenvalue \(\lambda=0\),
  and \(m\) be its root vector of order \(k\), i.e. \(a^{k-1}m\neq
  a^k m =0\). Then the restriction of \(\rho_a\) on the subspace
  generated by \(v_m\) is equivalent to the representation
  \(\rho_1^{k}\).
\end{prop}

\section{Spectrum and the Jordan Normal Form of a Matrix}
Now we are prepared to describe a spectrum of a matrix. Since the
functional calculus is an intertwining operator its support is
a decomposition into intertwining operators with prime
representations (we could not expect generally 
that these prime subrepresentations are irreducible).

Recall the transitive on \(\Space{D}{}\) group of inner
automorphisms of \(\SL\), which can send any
\(\lambda\in\Space{D}{}\) to \(0\) and are actually parametrised by
such a \(\lambda\).
This group extends Proposition~\ref{pr:Jordan-zero} to the complete
characterisation of \(\rho_a\) for matrices.
\begin{prop}  
  Representation \(\rho_a\) is equivalent to a direct sum of the
  prolongations \(\rho_1^{(k)}\) of \(\rho_1\) in the \(k\)th jet space
  \(\Space{J}{k}\) intertwined with inner automorphisms. Consequently
  the \textit{spectrum} of \(a\) (defined via the functional calculus
  \(\Phi=\oper{W}_m\)) labelled exactly by \(n\) pairs of numbers
  \((\lambda_i,k_i)\), \(\lambda_i\in\Space{D}{}\),
  \(k_i\in\Space[+]{Z}{}\), \(1\leq i \leq n\) some of whom could
  coincide.
\end{prop}
Obviously this spectral theory is a fancy restatement of the \emph{Jordan
  normal form} of matrices.

\begin{figure}[tb]
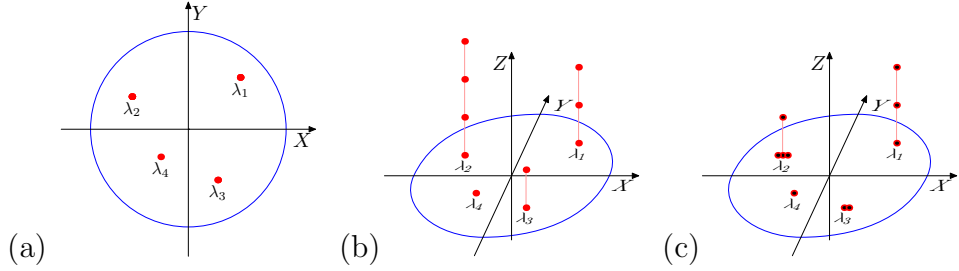

  \begin{center}
 (a) \includegraphics[scale=.65]{calc1vr.10}\hfill
  (b)\includegraphics[scale=.65]{calc1vr.27}\hfill
  (c)\includegraphics[scale=.65]{calc1vr.40}
 \caption[Three dimensional spectrum]{Classical spectrum (a) of a
   matrix vs. the new version (b) with its mapping (c).}
    \label{fig:3dspectrum}
  \end{center}
\end{figure}

\begin{example}
  \label{ex:3dspectrum}
  Let \(J_k(\lambda)\) denote the Jordan block of the length \(k\) for the
  eigenvalue \(\lambda\). On the Fig.~\ref{fig:3dspectrum} there are two
  pictures of the spectrum for the matrix
  \begin{displaymath}
    a=J_3\left(\lambda_1\right)\oplus     J_4\left(\lambda_2\right) 
    \oplus J_1\left(\lambda_3\right) \oplus      J_2\left(\lambda_4\right),
  \end{displaymath} 
  where
  \begin{displaymath}
    \lambda_1=\frac{3}{4}e^{i\pi/4}, \quad
    \lambda_2=\frac{2}{3}e^{i5\pi/6}, \quad
    \lambda_3=\frac{2}{5}e^{-i3\pi/4}, \quad
    \lambda_4=\frac{3}{5}e^{-i\pi/3}.
  \end{displaymath} Part (a) represents the conventional two-dimensional
  image of the spectrum, i.e. eigenvalues of \(a\), and
  \href{http://maths.leeds.ac.uk/~kisilv/calc1vr.gif}{(b) describes
  spectrum \(\spec{} a\) arising from the wavelet construction}. The
  first image did not allow to distinguish \(a\) from many other
  essentially different matrices, e.g. the diagonal matrix
  \begin{displaymath}
    \diag\left(\lambda_1,\lambda_2,\lambda_3,\lambda_4\right),
  \end{displaymath}
  which even have a different dimensionality.
  At the same time the Fig.~\ref{fig:3dspectrum}(b)
  completely characterise \(a\) up to a similarity. Note that each point of
  \(\spec a\) on Fig.~\ref{fig:3dspectrum}(b) corresponds to a particular
  root vector, which spans a primary subrepresentation.
\end{example}

\section{Spectral Mapping Theorem}
As was mentioned in the Introduction a resonable spectrum should be
linked to the corresponding functional calculus by an appropriate
spectral mapping theorem. The new version of spectrum is based on
prolongation of \(\rho_1\) into jet spaces (see
Section~\ref{sec:jet-bundl-prol-1}). Naturally a correct version of
spectral mapping theorem should also operate in jet spaces. 
%We give such a construction now.

Let \(\phi: \Space{D}{} \rightarrow \Space{D}{}\) be a holomorphic
map, let us define its action on functions \([\phi_*
f](z)=f(\phi(z))\). According to the general formula~\eqref{eq:prolong-def}
we can define the prolongation
\(\phi_*^{(n)}\) onto the jet space \(\Space{J}{n}\). Its associated
action \(\rho_1^k \phi_*^{(n)}=\phi_*^{(n)}\rho_1^n\) on the pairs
\((\lambda,k)\) is given by the formula:
\begin{equation}\label{eq:phi-star-action}
  \phi_*^{(n)}(\lambda,k)=\left(\phi(\lambda),
    \left[\frac{k}{\deg_\lambda \phi}\right]\right),
\end{equation}
where \(\deg_\lambda \phi\) denotes the degree of zero of the function
\(\phi(z)-\phi(\lambda)\) at the point \(z=\lambda\) and \([x]\) denotes
the integer part of \(x\). %We are ready to state

\begin{thm}[Spectral mapping] 
  Let \(\phi\) be a holomorphic mapping  \(\phi: \Space{D}{}
  \rightarrow \Space{D}{}\) and its prolonged action \(\phi_*^{(n)}\) defined
  by~\eqref{eq:phi-star-action}, then
  \begin{displaymath}
    \spec \phi(a) = \phi_*^{(n)} \spec a 
  \end{displaymath}
\end{thm} 

The explicit expression of~\eqref{eq:phi-star-action} for
\(\phi_*^{(n)}\), which involves derivatives of \(\phi\) upto \(n\)th order,
is known, see for example~\cite[Thm.~6.2.25]{HornJohnson94}, but was
not recognised before as form of spectral mapping.

\begin{example}
  Let us continue with Example~\ref{ex:3dspectrum}. Let \(\phi\) map
  all four eigenvalues \(\lambda_1\), \ldots, \(\lambda_4\) of the
  matrix \(a\) into themselves. Then Fig.~\ref{fig:3dspectrum}(a) will
  represent the classical spectrum of \(\phi(a)\) as well as \(a\).

  However Fig.~\ref{fig:3dspectrum}(c) shows mapping of the new
  spectrum for the case
  \(\phi\)  has
  \textit{orders of zeros} at these points as follows: the order \(1\)
  at \(\lambda_1\), exactly the order \(3\) at \(\lambda_2\), an order
  at least \(2\) at \(\lambda_3\), and finally any order at
  \(\lambda_4\).
\end{example}
\bibliographystyle{amsplain} 
\bibliography{abbrevmr,akisil,analyse,algebra}

\newcommand{\noopsort}[1]{} \newcommand{\printfirst}[2]{#1}
  \newcommand{\singleletter}[1]{#1} \newcommand{\switchargs}[2]{#2#1}
  \newcommand{\irm}{\textup{I}} \newcommand{\iirm}{\textup{II}}
  \newcommand{\vrm}{\textup{V}}
  \providecommand{\cprime}{'}\providecommand{\arXiv}[1]{\eprint{http://arXiv.o%
rg/abs/#1}{arXiv:#1}}
\providecommand{\bysame}{\leavevmode\hbox to3em{\hrulefill}\thinspace}
\providecommand{\MR}{\relax\ifhmode\unskip\space\fi MR }
\providecommand{\MRhref}[2]{%
  \href{http://www.ams.org/mathscinet-getitem?mr=#1}{#2}
}
\providecommand{\href}[2]{#2}
\begin{thebibliography}{10}

\bibitem{Anderson69}
Robert F.~V. Anderson, \emph{The {W}eyl functional calculus}, J. Functional
  Analysis \textbf{4} (1969), 240--267. \MR{58 \#30405}

\bibitem{HornJohnson94}
Roger~A. Horn and Charles~R. Johnson, \emph{Topics in matrix analysis},
  Cambridge University Press, Cambridge, 1994, Corrected reprint of the 1991
  original. \MR{95c:15001}

\bibitem{Kirillov76}
A.~A. Kirillov, \emph{Elements of the theory of representations},
  Springer-Verlag, Berlin, 1976, Translated from the Russian by Edwin Hewitt,
  Grundlehren der Mathematischen Wissenschaften, Band 220. \MR{54 \#447}

\bibitem{Kisil95i}
Vladimir~V. Kisil, \emph{M\"obius transformations and monogenic functional
  calculus}, \href{http://www.ams.org/era/}{Electron. Res. Announc. Amer. Math.
  Soc.} \textbf{2} (1996), no.~1,
  \href{http://www.ams.org/jourcgi/amsjournal?fn=120&pg1=pii&s1=S1079--6762--9%
6--00004--2}{26--33}, (electronic) \MR{98a:47018}.

\bibitem{Kisil97c}
\bysame, \emph{Analysis in{ $\Space{R}{1,1}$} or the principal function
  theory}, Complex Variables Theory Appl. \textbf{40} (1999), no.~2, 93--118,
  \arXiv{funct-an/9712003}. \MR{2000k:30078}.

\bibitem{Kisil97a}
\bysame, \emph{Two approaches to non-commutative geometry}, Complex Methods for
  Partial Differential Equations (H.~Begehr, O.~Celebi, and W.~Tutschke, eds.),
  Kluwer Academic Publishers, Netherlands, 1999, \arXiv{funct-an/9703001},
  \MR{2001a:01002}, pp.~219--248.

\bibitem{Kisil98a}
\bysame, \emph{Wavelets in {Banach} spaces}, Acta Appl. Math. \textbf{59}
  (1999), no.~1, 79--109, \arXiv{math/9807141}. \MR{2001c:43013}.

\bibitem{Kisil01a}
\bysame, \emph{Spaces of analytical functions and wavelets---{Lecture} notes},
  \arXiv{math.CV/0204018}, 2000--2002, 92 p.

\bibitem{Kisil02c}
\bysame, \emph{Meeting {Descartes} and {Klein} somewhere in a noncommutative
  space}, Procedsings of ICMP2000, 2002, \arXiv{math-ph/0112059}, p.~25~p.

\bibitem{Nikolskii86}
N.~K. Nikol{\cprime}ski{\u\i}, \emph{Treatise on the shift operator},
  Springer-Verlag, Berlin, 1986, Spectral function theory, With an appendix by
  S. V. Hru\v s\v cev [S. V. Khrushch\"ev] and V. V. Peller, Translated from
  the Russian by Jaak Peetre. \MR{87i:47042}

\bibitem{Olver93}
Peter~J. Olver, \emph{Applications of {L}ie groups to differential equations},
  second ed., Springer-Verlag, New York, 1993. \MR{94g:58260}

\bibitem{Olver95}
\bysame, \emph{Equivalence, invariants, and symmetry}, Cambridge University
  Press, Cambridge, 1995. \MR{96i:58005}

\end{thebibliography}

\end{document}